\documentstyle[12pt]{article}
\begin{document}
\bigskip
\noindent
{\bf REACTION-DIFFUSION SYSTEMS AND NONLINEAR WAVES}\par
\begin{center}
R.K. SAXENA\\
Department of Mathematics and Statistics, Jai Narain Vyas University\\
Jodhpur – 342004, India\\[0.5cm]
A.M. MATHAI\\
Department of Mathematics and Statistics, McGill University\\
Montreal, Canada H3A 2K6\\[0.5cm]
H.J. HAUBOLD\\
Office for Outer Space Affairs, United Nations\\
P.O.Box   500, A–1400 Vienna, Austria
\end{center}
\noindent
{\bf Abstract.}
The authors investigate the solution of a nonlinear reaction-diffusion equation connected with nonlinear waves. The equation discussed is more general than the one discussed recently by Manne, Hurd, and Kenkre (2000). The results are presented in a compact and elegant form in terms of Mittag-Leffler functions and generalized Mittag-Leffler functions, which are suitable for numerical computation. The importance of the derived results  lies in the fact that numerous results on fractional  reaction, fractional diffusion, anomalous diffusion problems,  and fractional  telegraph equations  scattered in the literature can be derived, as special cases,  of the  results investigated in this article.

\section{Introduction}

Reaction-diffusion models have found numerous applications in pattern formation in biology, chemistry, and physics, see Smoller (1983), Grindrod (1991), Gilding and  Kersner (2004), and Wilhelmsson and Lazzaro (2001). These systems show that diffusion can produce the spontaneous formation spatio-temporal patterns. For details, refer to the work of Nicolis and Prigogine (1977), and Haken (2004). A general model for reaction–diffusion systems is discussed by Henry and Wearne (2000, 2002), and Henry, Langlands, and Wearne (2005). A piecewise linear approach in connection with the diffusive processes has been developed by Strier, Zanette, and Wio (1995) which leads to analytic results in reaction-diffusion  systems. A similar approach was recently used by Manne, Hurd, and Kenkre (2000) to investigate effects on the propagation of nonlinear wave fronts.\\

The simplest reaction-diffusion models can be described by an equation

\begin{equation}
\frac{\partial{N}}{\partial{t}}=D\frac{\partial^2{N}}{\partial{x^2}}+\gamma F(N), 
\end{equation}
where $D$ is the diffusion coefficient and $F(N)$ is a nonlinear function representing reaction
kinetics. It is interesting to observe that for  $F(N) = \gamma N(1-N),$ eq.(1) reduces to Fisher-Kolmogorov equation and if we set  $F(N) = \gamma N(1-N^2)$, it reduces to  the real Ginsburg-Landau equation. \\

A generalization of (1) has been considered by Manne, Hurd, and Kenkre (2000) in the form

\begin{equation}
\frac{\partial^2{N}}{\partial{t^2}}+a\frac{\partial{N}}{\partial{t}}=\nu^2\frac{\partial^2{N}}{\partial{x^2}}+\xi^2 N(x,t),
\end{equation}
where $\xi$ indicates the strength of the nonlinearity of the system.
In this article, we present a straightforward method for the systematic derivation of the solution of nonlinear reaction-diffusion equations connected with nonlinear waves, which is more general than the equation (2). The results are derived in a closed-form, by the application of Laplace and Fourier transforms, which are suitable for numerical computation. The present study is in continuation of our investigations reported earlier in the articles (Saxena, Mathai, and Haubold, 2002, 2004, 2004a, 2004b, 2005).

\section{ Mathematical Prerequisites} 
A generalization of the Mittag-Leffler function (Mittag-Leffler 1903, 1905)
\begin{equation}
E_\alpha(z)=\sum^\infty_{n=0}\frac{z^n}{\Gamma(n\alpha+1)},\;\;\alpha \in C, Re
(\alpha)>0,
\end{equation}                                                     
was introduced by Wiman (1905) in the general form 
\begin{equation}
E_{\alpha,\beta}(z)=\sum^\infty_{n=0}\frac{z^n}{\Gamma(n\alpha+\beta)},\;\alpha,\beta\in C, Re(\alpha)>0, Re(\beta)>0.
\end{equation}

The main results of these functions are available in the handbook of Erd\'{e}lyi, Magnus, Oberhettinger, and Tricomi (1955, Section 18.1) and the monographs  by Dzherbashyan (1966, 1993). Prabhakar (1971) introduced a generalization of (4) in the form 

\begin{equation}
E^\gamma_{\alpha,\beta}(z)=\sum^\infty_{n=0}\frac{(\gamma)_nz^n}{\Gamma(n\alpha+\beta)(n)!}\;\alpha,\beta,\gamma \in C; Re(\alpha), Re(\beta), Re(\gamma)>0,
\end{equation}                                                             
where $(\gamma)_r$ is Pochhammer's  symbol, defined by 
\begin{equation}
(\gamma)_0=1, (\gamma)_r=\gamma(\gamma+1)(\gamma+2)\ldots(\gamma+r-1), r=1,2,\ldots, \gamma\neq 0.
\end{equation} 
It is an entire function with $\rho = [Re(\alpha)]^{-1}$    (Prabhakar, 1971). The solution of  generalized Volterra-type differ-integral equations associated with this function as a kernel is derived by Kilbas, Saigo and Saxena 
(2004). A general theory of generalized fractional calculus based on this function has been developed by Kilbas, Saigo and Saxena (2004), generalizing the results for Riemann-Liouville fractional integrals and derivatives, which form the backbone of fractional differ-integral equations.

For $\gamma=1$, this function coincides with (4), while for $\beta=\gamma=1$ with (3)
\begin{equation}
E^1_{\alpha,\beta} = E_{\alpha,\beta}(z), E^1_{\alpha,1}(z) = E_\alpha(z). 
\end{equation}
We also have 
\begin{equation}
\Phi(\alpha,\beta; z)= \;_1F_1(\alpha,\beta;z)= \Gamma(\beta)E_{1,\beta}^\alpha (z),
\end{equation}                                                       
where $\Phi(\alpha,\beta;z)$  is Kummer's confluent hypergeometric function defined in Erd\'{e}lyi, Magnus, Oberhettinger, and Tricomi (1953, p.248, eq.(1)). Prabhakar (1971, p.8, eq.(2.5)) has shown that
\begin{equation}
\int_0^\infty t^{\gamma-1}e^{-st}E^\delta_{\beta,\gamma}(\omega t^\beta)dt = s^{-\gamma}(1-\omega s^{-\beta})^{-\delta},
\end{equation}
where $Re(\beta)>0, Re(\gamma)>0, Re(s)> 0\;\mbox{and}\; s>|\omega|^{\frac{1}{Re(\beta)}}.$

The Riemann-Liouville fractional integral of order $\nu$ is defined by 
Miller and Ross (1993, p.45) 
\begin{equation}
_0D_t^{-\nu}f(t)=\frac{1}{\Gamma(\nu)}\int_0^t(t-u)^{\nu-1}f(u)du, Re(\nu)>0.
\end{equation}

Following Samko, Kilbas, and Marichev (1990, p.37), we define the fractional derivative for $Re(\alpha)>0$ in the form 
\begin{equation}
_0D_t^\alpha f(t)=\frac{1}{\Gamma(n-\alpha)}\frac{d^n}{dt^n}\int_0^t\frac{f(u)du}{(t-u)^{\alpha-n+1}}\,;n=[\alpha]+1,
\end{equation}
where $[\alpha]$ means the integral part of  the number $\alpha$.
In particular, if $\;0 <\alpha < 1,$
\begin{equation}
_0D_t^\alpha f(t)= \frac{d}{dt}\frac{1}{\Gamma(1-\alpha)}\int_0^t\frac{f(u)du}{(t-u)^\alpha},
\end{equation}
and if $\alpha = n \in N = \left\{1,2,\ldots \right\},$ then
\begin{equation}
_0D_t^n f(t)=D^nf(t),(D=d/dt),
\end{equation}
which is the standard derivative of order n.
From Erd\'{e}lyi, Magnus, Oberhettinger, and Tricomi 
(1954a, 1954b, p.182), we have
\begin{equation} 
L\left\{_0D_t^{-\nu} f(t); s\right\}= s^{-\nu}\tilde{f}(s),
\end{equation}                                                                                 
where $\tilde{f}(s)$  is the Laplace transform of  $f(t)$, defined by 
\begin{equation}
L\left\{f(t);s\right\}=\tilde{f}(s) =\int^\infty_0 \mbox{exp}(-st) f(t)dt, Re(s)>0,
\end{equation}                                                                             
which may be written symbolically, as follows
\begin{equation}
\tilde{f}(s)=L\left\{f(t);s\right\}\mbox{or}\;\;  f(t)= L^{-1}\left\{\tilde{f}(s);t\right\},
\end{equation}                                                                    
provided that the function f(t) is continuous for $t\geq 0$ and of exponential order as $t\rightarrow \infty$. The Laplace transform of the fractional derivative is given by  
Oldham and Spanier (1974, p.134, Eq. (8.1.3))
\begin{equation}
L\left\{_0D_t^\alpha f(t);s\right\}=s^\alpha \tilde{f}(s)+\sum^n_{r=1}s^{r-1}\;_0D_t^{\alpha-r} f(t)|_{t=0}.
\end{equation}                                                   
In certain boundary-value problems,  the following fractional derivative of order $\alpha > 0$ of a causal function  $f(t)$  (that is, $f(t) = 0,$ for $t < 0$ ), is introduced by Caputo (1969) in the form
\begin{eqnarray}
D_t^\alpha f(t)&=&\frac{1}{\Gamma(m-\alpha)}\int^t_0\frac{f^{(m)}(\tau)d\tau}{(t-\tau)^{\alpha-m+1}},\;m-1<\alpha\leq m, Re(\alpha)>0, m\in N,.\nonumber\\
&=&\frac{d^mf}{dt^m}, \;\mbox{if}\; \;\alpha =m,
\end{eqnarray}
where $\frac{d^m}{dx^m} f$ is the $m^{th}$ derivative of $f$.
                                                                         
Caputo (1969) has given the Laplace transform of the fractional derivative as 
\begin{equation}
L \left\{D_t^\alpha f(t); s\right\}  =s^\alpha F(s)-\sum^{m-1}_{r=0}s^{\alpha-r-1}f^{(r)}(0+), m-1<\alpha\leq m,
\end{equation}
where $F(s)$ is the Laplace transform of  $f(t)$.\\

The above formula is useful in deriving the solution of differ-integral equations of fractional order governing certain physical problems of reaction and diffusion. 
We also need the Weyl fractional operator defined by
\begin{equation}
_{-\infty}D_x^\mu f(t)= \frac{1}{\Gamma(n-\mu)}\frac{d^n}{dt^n}\int_{-\infty}^t\frac{f(u)du}{(t-u)^{\mu-n+1}},
\end{equation}
where $n = [\mu]$  is an integral  part of $\mu>0.$
Its Fourier transform is given by Metzler and Klafter (2000, p.59, A.11, 2004)
\begin{equation}
F\left\{_{-\infty}D_x^\mu f(x)\right\}=(ik)^\mu f^*(k),
\end{equation}                                                                 
where we define the  Fourier transform by the integral equation  
\begin{equation}
h^*(q)=\int_{-\infty}^\infty h(x)exp(iqx)dx.
\end{equation}
Following the convention initiated by Compte (1996), we suppress the imaginary unit in Fourier space by adopting a slightly modified form of the result (21) in our investigations (Metzler and Klafter 2000, p. 59, A.12, 2004)
\begin{equation}
F\left\{_{-\infty}D_x^\mu f(x)\right\}=-|k|^\mu f^*(k).
\end{equation}      
Now we will establish the following results, which provide the inverse Laplace transforms of certain algebraic functions and are directly applicable in the analysis of reaction-diffusion systems that follows.\par
It will be shown here that
\begin{equation}
(A) \;L^{-1}\left\{\frac{s^{\alpha-1}}{s^\alpha+as^\beta+b};t\right\}=\sum^\infty_{r=0}(-a)^r t^{(\alpha-\beta)r}E^{r+1}_{\alpha,(\alpha-\beta)r+1}(-bt^\alpha),
\end{equation}
where $Re(\alpha)>0, Re(\beta)>0,Re(s)>0,\;\;|\frac{as^\beta}{s^\alpha+b}|<1$, and $E^\delta_{\beta,\gamma}(x)$ is the generalized Mittag-Leffler function defined in eq. (5).\par
\noindent
{\bf Proof.} We have 
      
\begin{eqnarray} 
\frac{s^{\alpha-1}}{s^\alpha+as^\beta+b}&=&\frac{s^{\alpha-\beta-1}}{(s^{\alpha-\beta}+bs^{-\beta})(1+\frac{a}{s^{\alpha-\beta}+bs^{-\beta}})}\nonumber\\                                        
&=&s^{\alpha-\beta-1}\sum^\infty_{r=0}\frac{(-a)^r}{(s^{\alpha-\beta}+bs^{-\beta})^{r+1}}\nonumber\\
&=&\sum^\infty_{r=0}\frac{(-a)^r s^{\alpha+\beta(r+1)-\beta-1}}{(s^\alpha+b)^{r+1}}.
\end{eqnarray}
Taking the inverse Laplace transform  of the result (25) and using eq. (9), we obtain 
the result (24). The term by term inversion is justified by virtue of a theorem given by Doetsch (1956, § 22). 
In a similar manner, the following results can be established
\begin{eqnarray}
(B) &L^{-1}&\left\{\frac{s^{\beta-1}}{s^\alpha+as^\beta+b};t\right\}\nonumber\\
&=& t^{\alpha-\beta}\sum^\infty_{r=0}(-a)^r t^{(\alpha-\beta)r}E^{r+1}_{\alpha,(\alpha-\beta)(r+1)+1}(-bt^\alpha),
\end{eqnarray}
where $Re(\alpha)>0, Re(\beta)>0, Re(s)>0,|\frac{as^\beta}{s^\alpha+b}|<1,$  and $\alpha >\beta$.
\begin{eqnarray}
(C) &L^{-1}&\left\{\frac{1}{s^\alpha+as^\beta+b};t\right\}\nonumber\\
&=&t^{\alpha-1}\sum^\infty_{r=0}(-a)^rt^{(\alpha-\beta)r}E^{r+1}_{\alpha, \alpha+(\alpha-\beta)r} (-bt^\alpha),
\end{eqnarray}      
where  
$Re(\alpha)>0, Re(\beta)>0, Re(s)>0, and \;|\frac{as^\beta}{s^\alpha+b}|<1$.\\
We now show that
\begin{eqnarray}
(D) &L^{-1}&\left\{\frac{s^{2\alpha-1}+as^{\alpha-1}}{s^{2\alpha}+as^\alpha+b};t\right\}\nonumber\\
&=&\frac{1}{\sqrt{(a^2-4b)}}[(\lambda+a)E_\alpha(\lambda t^\alpha)-(\mu+a)E_\alpha(\mu t^\alpha)],
\end{eqnarray}  
where $a^2-4b>0$ and $E_\alpha(x)$  is the Mittag-Leffler function defined in eq. (3), $Re(\alpha)>0, Re(s)>0, $ and  $\lambda$ and 
$\mu$ are the real and distinct  roots of the quadratic equation,
$x^2+ax+b=0,$
namely $\lambda=\frac{1}{2}(-a+\sqrt{(a^2-4b)})$ and $\mu=\frac{1}{2}(-a-\sqrt{(a^2-4b)})$.\par
\noindent
{\bf Proof.} We have 
\begin{eqnarray}
\frac{s^{2\alpha-1}+as^{\alpha-1}}{s^{2\alpha}+as^\alpha+b}&=&\frac{1}{\lambda-\mu}[\frac{s^{2\alpha-1}+as^{\alpha-1}}{s^\alpha-\lambda}-\frac{s^{2\alpha-1}+as^{\alpha-1}}{s^\alpha-\mu}]\nonumber\\
&=& \frac{1}{\lambda-\mu}[\frac{(\lambda+a)s^{\alpha-1}}{s^\alpha-\lambda}-\frac{(\mu+a)s^{\alpha-1}}{s^\alpha-\mu}].
\end{eqnarray}     
Taking the inverse Laplace transform of (29) gives 
\begin{eqnarray*}
&L^{-1}&\left\{\frac{s^{2\alpha-1}+as^{\alpha-1}}{s^{2\alpha}+as^\alpha+b}\right\}\nonumber\\
&=&\frac{1}{\lambda-\mu}[(\lambda+a)E_\alpha(\lambda t^\alpha)-(\mu+a)E_\alpha(\mu t^\alpha)].
\end{eqnarray*}   
This completes the proof of eq. (28).\\ 
In a similar manner, it can be shown that 
\begin{equation}
(E)\; L^{-1}\left\{\frac{1}{s^{2\alpha}+as^\alpha+b}\right\}=\frac{t^{\alpha-1}}{\lambda-\mu}[E_{\alpha,\alpha}(\lambda t^\alpha)-E_{\alpha,\alpha}(\mu t^\alpha)],
\end{equation}
where $Re(\alpha) > 0, Re(s) > 0$, and $\lambda$ and $\mu$ are given along with (28).

\section{ Solution of Fractional Reaction-Diffusion Equation}

In this section, it is proposed to derive the solution of the fractional-diffusion system connected with nonlinear waves governed by the eq. (31). This system is a generalized form of the reaction-diffusion equation recently studied by Manne, Hurd, and Kenkre (2000). The result is given in the form of the following theorem.\\
\noindent
{\bf Theorem.} Consider the fractional reaction-diffusion equation 
\begin{equation} 
_0D_t^\alpha N(x,t)+a\;_0D_t^\beta N(x,t)
= \nu^2 _{-\infty}D_x^\gamma N(x,t)+\xi^2 N(x,t)+\varphi(x,t),
\end{equation}
$$x\in \Re,t>0,0\leq\alpha\leq 1,0\leq \beta \leq 1,$$
with initial conditions 
\begin{equation}
N(x,0) = f(x),\; \mbox{for} \;(x\in \Re),
\end{equation}                                                                                 
where $\nu^2$  is a diffusion coefficient, $\varphi$  is a constant which describes the nonlinearity in the system, and $\varphi(x,t)$ is a nonlinear function for reaction kinetics, then there holds the following formula for the solution of (31)
\begin{eqnarray}
N(x,t)&=&
\sum^\infty_{r=0}\frac{(-a)^r}{2\pi}\int^\infty_{-\infty}t^{(\alpha-\beta)r}\;f^*(k)exp(-kx)\times\nonumber\\
&\times& [E^{r+1}_{\alpha,(\alpha-\beta)r+1}(-bt^\alpha)+t^{\alpha-
\beta}
E^{r+1}_{\alpha,(\alpha-\beta)(r+1)+1}(-bt^\alpha)]dk\nonumber\\
&+& \sum^\infty_{r=0}\frac{(-a)^r}{2\pi}\int^t_0\xi^{\alpha+(\alpha-\beta)r-1}\int^\infty_{-\infty}\varphi(k,t-\xi)exp(-ikx)\times\nonumber\\
&\times&E^{r+1}_{\alpha, \alpha+(\alpha-\beta)r}(-b\xi^\alpha)dkd\xi,
\end{eqnarray}
where $\alpha>\beta$  and $E^\delta_{\beta,\gamma}(.)$  is the generalized Mittag-Leffler function, defined in (5) and  
$b=\nu^2|k|^\gamma-\xi^2$.\\
\noindent
{\bf Proof.} Applying the Laplace transform with respect to the time variable
 t and using the boundary conditions , we find that 
\begin{eqnarray}
&&s^\alpha \tilde{N}(x,s)-s^{\alpha-1}f(x)+as^\beta \tilde{N}(x,s)-as^{\beta-1}f(x)\nonumber\\
&=&\nu^2\;_{-\infty}D_x^\nu \tilde{N}(x,s)+\xi^2 \tilde{N}(x,s)+ \tilde{f}(x,s).
\end{eqnarray}
If we apply the Fourier transform with respect to the space variable $x$, it yields           
\begin{eqnarray}
&&s^\alpha \tilde{N}^*(k,s)-s^{\alpha-1}f^*(k)+as^\beta \tilde{N}^*(k,s)-as^{\beta-1}f^*(k)\nonumber\\
&=&-\nu^2|k|^\gamma\tilde{N^*}(k,s)+\xi^2\tilde{N^*}(k,s)+\tilde{f^*}(k,s).
\end{eqnarray}       
                                             
Solving for $\tilde{N^*}(k,s)$, it gives
\begin{equation}
\tilde{N^*}(k,s)=\frac{(s^{\alpha-1}+as^{\beta-1})f^*(k)+\tilde{f^*}(k,s)}{s^\alpha+as^\beta+b},
\end{equation}
where $b= \nu^2|k|^\gamma-\xi^2$.
To invert eq. (36), it is convenient to first invert the Laplace transform and then the Fourier transform. Inverting the Laplace transform with the help of the results (28) and (30), yields
\begin{eqnarray}
N^*(k,t)&=&\sum^\infty_{r=0}(-a)^r t^{(\alpha-\beta)r} f^*(k)[E^{r+1}_{\alpha,(\alpha-\beta)r+1}(-bt^\alpha)+t^{\alpha-\beta}\times\nonumber\\
&\times&E^{r+1}_{\alpha,(\alpha-\beta)(r+1)+1}(-bt^\alpha)]\nonumber\\
&+&\sum^\infty_{r=0}(-a)^r\int_0^t\varphi^*(k,t-\xi)\xi^{\alpha+(\alpha-\beta)r-1}\times\nonumber\\
&\times& E^{r+1}_{\alpha,(\alpha-\beta)r+\alpha}(-b\xi^\alpha)d\xi.
\end{eqnarray}
                                                                                        
Finally, the inverse Fourier transform gives the desired solution in the form 
\begin{eqnarray}
N(x,t)&=&\sum^\infty_{r=0}\frac{(-a)^r}{2\pi}\int^\infty_{-\infty}t^{(\alpha-\beta)r}f^*(k)[E^{r+1}_{\alpha(\alpha-\beta)r+1}(-bt^\alpha)+t^{\alpha-\beta}\times\nonumber\\
&\times&E^{r+1}_{\alpha(\alpha-\beta)(r+1)+1}(-bt^\alpha)] exp(-ikx)dk\nonumber\\
&+&\sum^\infty_{r=0}\frac{(-a)^r}{2\pi}\int_0^t\xi^{\alpha+(\alpha-\beta)r-1}\int^\infty_{-\infty}exp(-ikx)\varphi^*(k,t-\xi)\times\nonumber\\
&\times& E^{r+1}_{\alpha,\alpha+(\alpha-\beta)r}(-b\xi^\alpha)dkd\xi.
\end{eqnarray}
This completes the proof of the theorem.

\section{Special cases}

When $f (x) =\delta(x)$, where $\delta(x)$ is the Dirac delta function, the theorem reduces to the following\\
\noindent 
{\bf Corollary 1.} Consider the fractional reaction-diffusion system 
\begin{equation}
_0D_t^\alpha N(x,t)+a\;_0D_t^\beta N(x,t)
= \nu^2\;_{-\infty}D^\gamma_x N(x,t)+\xi^2 N(x,t)+\varphi(x,t),
\end{equation}
subject to the initial conditions 
\begin{equation}
N(x,0) = \delta(x)\; \mbox{for} \;0\leq a \leq 1\; \mbox{and}\;\;0\leq\beta\leq 1,
\end{equation}
where $\delta(x)$   is the Dirac delta function. Here $\xi$   is a constant that describes the nonlinearity in the system, and $\varphi(x,t)$  is a  nonlinear function which belongs to the reaction kinetics. Then there exists the following eq. for the solution of (39), subject to the initial conditions (40)
\begin{eqnarray}
N(x,t)&=&\sum^\infty_{r=0}\frac{(-a)^r}{2\pi}\int^\infty_{-\infty} e^{-ikx}[t^{(\alpha-\beta)r}E^{r+1}_{\alpha,(\alpha-\beta)r+1}(-bt^\alpha)+t^{(\alpha-\beta)(r+1)}\times\nonumber\\
&\times&E_{\alpha(\alpha-\beta)(r+1)+1}(-bt^\alpha)]dk\nonumber\\
&+&\sum^\infty_{r=0}\frac{(-a)^r}{2\pi}\int_0^t\xi^{\alpha+((\alpha-\beta)r-1}\int^\infty_{-\infty}\varphi^*(k,t-\xi)\nonumber\\
&\times& exp(-ikx)E^{r+1}_{\alpha,(\alpha-\beta)r+\alpha}(-b\xi^\alpha)dkd\xi,
\end{eqnarray}
where $b =  \nu^2|k|^\gamma -\xi^2$.

Now if we set $f(x) =\delta(x), \gamma=2, \alpha$ is replaced by $2\alpha$  and $\beta$ by $\alpha$, the following result is obtained.\\
\noindent
{\bf Corollary 2.} Consider the following reaction-diffusion system
$$\frac{\partial^{2\alpha} N(x,t)}{\partial t^{2\alpha}}+a\frac{\partial^\alpha N(x,t)}{\partial t^\alpha}=\nu^2\frac{\partial^2 N(x,t)}{\partial x^2}+\xi^2N(x,t)+\varphi(x,t)$$       
with the initial conditions
\begin{equation}
N(x,0)=\delta(x), \; N_t(x,0)=0,\;\;\; 0\leq\alpha\leq 1,
\end{equation}
$\varphi(x,t)$  is a nonlinear function belonging to the reaction kinetics. Then for the solution of (39) subject to the initial conditions (40), there holds the formula 
\begin{eqnarray}
N(x,t)&=&\frac{1}{2\pi\sqrt{(a^2-4b)}}\left[\int^\infty_{-\infty}exp(-ikx)\times \right.\\
&\times&\left\{(\lambda+a)E_\alpha(\lambda t^\alpha)-(\mu+a)E_\alpha(\mu t^\alpha)\right\}dk\nonumber\\
&+&\frac{1}{2\pi}\int_0^t\xi^{\alpha-1}\int^\infty_{-\infty}exp(-ikx)\varphi^*(k,t-\xi)\times\nonumber\\
&\times&\left.[E_{\alpha,\alpha}(\mu\xi^\alpha)-E_{\alpha, \alpha}(\mu\xi^\alpha)]dkd\xi\right],\nonumber
\end{eqnarray}                                                                                                                     
where $\lambda$  and  $\mu$ are the real and distinct roots of  the quadratic equation 
\begin{equation}
y^2+ay+b=0,
\end{equation}
given by 
\begin{equation}
\lambda=\frac{1}{2}\left(-a+\sqrt{(a^2-4b)}\right)\;\mbox{and}\; \mu=\frac{1}{2}\left(-a-\sqrt{(a^2-4b)}\right),
\end{equation} 
 where $b^2=\nu^2k^2-\xi^2$.\\
\noindent
{\bf Proof.} In order to prove (43), we replace $\alpha$  by $2\alpha$ and $\beta$  by $\alpha$, then eq. (36) becomes 
\begin{equation}
\tilde{N^*}(k,s)= \frac{s^{2\alpha-1}+as^{\alpha-1}+\tilde{\varphi^*}(k,s)}{s^{2\alpha}+as^\alpha+b}.
\end{equation}
                                                            
Taking the inverse Laplace transform and using the results (26) and
(30), yields 

\begin{eqnarray}
N^*(k,t)&=&\frac{1}{\lambda-\mu}[(\lambda+a)E_\alpha(\lambda t^\alpha)-(\mu+a)E_\alpha(\mu t^\alpha)]\\
&+&\int_0^t\varphi^*(k,t-\xi)\xi^{\alpha-1}[E_{\alpha, \alpha}(\lambda \xi^\alpha)-E_{\alpha, \alpha}(\mu\xi^\alpha)]d\xi,\;\;\lambda\neq \mu,\nonumber
\end{eqnarray}
where $\lambda$   and $\mu$ are given  in (45).                                                                                                        
The application of the inverse Fourier transform to the above equation
gives the desired result (43). 

Next, if we set $\varphi(x,t)=0,$ $\gamma=2$, replace $\alpha$ by $2\alpha$, and $\beta$  by $\alpha$  in\\  
(31), we then obtain the following result , which includes many known results on the fractional telegraph equations including the one recently given by Orsingher and Beghin (2004).\\
\noindent
{\bf Corollary 3.} Consider the following reaction-diffusion system

\begin{equation}
\frac{\partial^{2\alpha}N(x,t)}{\partial t^{2\alpha}}+a\frac{\partial^\alpha N(x,t)}{\partial t^\alpha}=\nu^2\frac{\partial^2 N(x,t)}{\partial x^2}+\xi^2 N(x,t),
\end{equation}
with the initial conditions
\begin{equation}
N(x,0)=\delta(x),\;\;\;N_t(x,0)=0,\;\;\;0\leq\alpha\leq1.
\end{equation} 
Then for the solution of (48), subject to the initial conditions (49), there holds the formula
\begin{eqnarray}
N(x,t)&=&\frac{1}{2\pi\sqrt{(a^2-4b)}}\times\\
&\times&\left[\int^{+\infty}_{-\infty} exp(-ikx)\left\{(\lambda+a)E_\alpha(\lambda t^\alpha)-(\mu+a)E_\alpha(\mu t^\alpha)\right\}dk\right],\nonumber
\end{eqnarray}
where $\lambda$  and $\mu$ are defined in (45), $b=\nu^2k^2-\xi^2$ and $E_\alpha(x)$ is the Mittag-Leffler function defined by (3). 
If we set $\xi^2=0$, then  corollary 3 reduces to the  result, which states that
the reaction-diffusion system 
\begin{equation}
\frac{\partial^{2\alpha}N(x,t)}{\partial t^{2\alpha}}+a\frac{\partial^\alpha N(x,t)}{\partial t^\alpha}=\nu^2\frac{\partial^2N(x,t)}{\partial x^2}
\end{equation}                       
with the initial conditions
\begin{equation}
N(x,0)=\delta(x),\;\;N_t(x,0)=0,\;\;0\leq \alpha \leq 1,
\end{equation}
has the solution, given by  
\begin{eqnarray}
N(x,t)&=&\frac{1}{2\pi\sqrt{(a^2-4b)}}\times\\
&\times&\left[\int^\infty_{-\infty} exp(-ikx)\left\{(\lambda+a)E_\alpha(\lambda t^\alpha)-(\mu+a)E_\alpha(\mu t^\alpha)\right\}dk\right],\nonumber
\end{eqnarray}
where $\lambda$   and $\mu$   are defined in (45), $b=\nu^2k^2$ and $E_\alpha(x)$ is the Mittag-Leffler function defined by (3).
Eq. (53) can be rewritten in the form 
\begin{eqnarray}
N(x,t)&=&\frac{1}{4\pi}\left[\int^\infty_{-\infty}exp(-ikx)\left\{(1+\frac{a}{\sqrt{(a^2-4\nu^2k^2)}})\right.\right.\\
&\times& \left.\left. E_\alpha(\lambda t^\alpha)+(1-\frac{a}{\sqrt{(a^2-4\nu^2k^2)}})E_\alpha(
\mu t^\alpha)\right\}dk\right],\nonumber
\end{eqnarray}
where $\lambda$  and $\mu$  are defined in (45) and $E_\alpha(x)$  is the Mittag-Leffler function, defined by (3). 
The equation (54) represents the solution of the time-fractional telegraph equation (51), subject to the initial conditions (52), recently solved by Orsingher and Beghin (2004). It may be remarked here that the solution as given by Orsingher and Beghin (2004) is in terms of the Fourier transform of the solution in the form given below.
The Fourier transform of the solution of the equations (51) and 
(52)  can be expressed in the form 
\begin{equation}
N^*(x,t)=\frac{1}{2}\left\{(1+\frac{a}{\sqrt{(a^2-4\nu^2k^2})})E_\alpha (\lambda t^\alpha)+(1-\frac{a}{\sqrt{(a^2-4\nu^2k^2))}})E_\alpha(\mu t^\alpha)\right\},
\end{equation}
where $\lambda$   and  $\mu$  are defined in (45) and $E_\alpha(x)$  is the Mittag-Leffler function defined by (3).
Finally, it is interesting to observe that the solution of  various  fractional-reaction and fractional diffusion  and fractional telegraph equations scattered in the literature can be derived as special cases of the theorem  established in this article. 
\section{Conclusions}
There is a host of reaction-diffusion equations such as eqs. (1) and (2) that allow the formation of wave fronts which maintain their shape despite the diffusive element of evolution contrary to linear expectation. Such generation of waves plays a particularly important role in spatio-temporal processes in physical systems, including astrophysical fusion plasmas (Kulsrud, 2005; Wilhelmsson and Lazarro, 2001; Gilding and Kersner 2004). Much work has been done for the numerical treatment of such equations. 

The motivation for the research reported in this paper is the derivation of analytic closed-form solutions of fractional reaction-diffusion equations (31), (39), (42), and (48) that give rise to nonlinear waves in a respective physical medium. For this purpose the paper summarizes specific techniques for Laplace, Fourier, and Mellin transforms, results for Mittag-Leffler functions, as well as the applications of Riemann-Liouville, Weyl, and Caputo fractional calculus for tackling fractional reaction-diffusion equations. Closed-form solutions of the equations are given in terms of Fox's function and their behavior for small and large values of the respective parameter is derived. 
\bigskip
\noindent
\section*{References}\par
\bigskip
\noindent
Caputo, M.: 1969, \emph {Elasticita e Dissipazione}, Zanichelli, Bologna.\\
Doetsch,G.: 1956, \emph {Anleitung zum Praktischen Gebrauch}

\emph {der Laplace-Transformation}, Oldenbourg, Munich. \\
Dzherbashyan, M.M.: 1966, \emph {Integral Transforms and Representation of}

\emph {Functions in Complex Domain} (in Russian), Nauka, Moscow.\\
Dzherbashyan, M.M.: 1993, \emph {Harmonic Analysis and Boundary Value}

\emph {Problems in the Complex Domain}, Birkhaeuser-Verlag, Basel.\\
Erd\'{e}lyi, A., Magnus, W., Oberhettinger, F., and Tricomi, F.G.: 1953,

\emph {Higher Transcendental Functions}, Vol. {\bf 1}, McGraw-Hill, New York, 

Toronto, and London.\\
Erd\'{e}lyi, A., Magnus, W., Oberhettinger, F., and Tricomi, F.G.: 1954a, \emph {Tables}

\emph {of Integral Transforms}, Vol. {\bf 1}, McGraw-Hill, New York, Toronto,

and London.\\
Erd\'{e}lyi, A., Magnus, W., Oberhettinger, F., and Tricomi, F.G.: 1954b,

\emph {Tables of Integral Transforms}, Vol. {\bf 2}, McGraw-Hill, New York, Toronto,

and London.\\
Erd\'{e}lyi, A., Magnus, W., Oberhettinger, F., and Tricomi, F.G.: 1955,

\emph {Higher Transcendental Functions}, Vol. {\bf 3}, McGraw-Hill, New York, Toronto,

and London.\\
Gilding, B.H. and Kersner, R.: 2004, \emph {Travelling Waves in Nonlinear} 

\emph {Diffusion-Convection Reaction}, Birkhaeuser-Verlag, Basel-Boston- 

Berlin.\\
Grindrod, P.: 1991, \emph {Patterns and Waves: The Theory and Applications} 

\emph {of Reaction-Diffusion Equations}, Clarendon Press, Oxford.\\
Haken, H.: 2004, \emph {Synergetics: Introduction and Advanced Topics}, 

Springer-Verlag, Berlin-Heidelberg.\\
Henry, B.I and Wearne, S.L.: 2000, Fractional reaction-diffusion,

\emph {Physica A} {\bf 276}, 448-455.\\
Henry, B.I. and Wearne, S.L.: 2002, Existence of Turing instabilities 

in a two-species fractional reaction-diffusion system, 

\emph {SIAM Journal of Applied Mathematics} {\bf 62}, 870-887.\\
Henry, B.I., Langlands, T.A.M., and Wearne, S.L.: 2005, Turing 

pattern formation in fractional activator-inhibitor systems, 

\emph {Physical Review E} {\bf 72}, 026101.\\
Kilbas, A.A. and Saigo, M.:  2004, \emph {H-Transforms: Theory and}

\emph {Applications}, Chapman and Hall/CRC,  New York.\\
Kilbas, A.A., Saigo, M., and Saxena, R.K.: 2004, Generalized Mittag-
 
Leffler function  and generalized fractional calculus, \emph {Integral Transforms} 

\emph {and Special Functions} {\bf 15}, 31-49.\\
Kulsrud, R.M.: 2005, \emph {Plasma Physics for Astrophysics}, Princeton 

University Press, Princeton and Oxford.\\
Manne, K.K., Hurd, A.J., and Kenkre, V.M.: 2000, Nonlinear waves in

reaction-diffusion systems: The effect of transport memory,

\emph {Physical Review E} {\bf 61}, 4177-4184.\\
Metzler, R. and Klafter, J.: 2000, The random walk's guide to anomalous

diffusion: A fractional dynamics approach, \emph {Physics Reports} {\bf 339}, 1-77.\\
Metzler, R. and Klafter, J.: 2004, The restaurant at the end of the 

random walk: Recent developments in the description of anomalous 

transport by fractional dynamics, \emph {Journal of Physics A: Math. Gen.} {\bf 37}, 

R161-R208.\\
Miller, K.S. and Ross, B.: 1993, \emph {An Introduction to the Fractional Calculus}

\emph {and Fractional Differential Equations}, John Wiley and Sons, New York.\\
Mittag-Leffler, M.G.: 1903, Sur la nouvelle fonction $E_\alpha(x)$,

\emph {Comptes Rendus Acad.Sci. Paris (Ser.II)} {\bf 137}, 554-558.\\
Mittag-Leffler, M.G.: 1905, Sur la representation analytique d'une 

branche uniforme d'une fonction monogene,

\emph {Acta Mathematica} {\bf 29}, 101-181.\\
Nicolis, G. and Prigogine, I.: 1977, \emph {Self-Organization in Nonequilibrium} 

\emph {Systems: From Dissipative Structures to Order Through Fluctuations}, 

John Wiley and Sons, New York.\\  
Oldham, K.B. and Spanier, J.: 1974, \emph {The Fractional Calculus:}

\emph {Theory and Applications of Differentiation and Integration to Arbitrary}

\emph {Order}, Academic Press, New York; and Dover Publications, New York 

2006.\\
Orsingher, E. and Beghin, L.: 2004, Time-fractional telegraph equations

and telegraph processes with Brownian time, \emph {Probability Theory}

\emph {and Related Fields} {\bf 128}, 141-160.\\
Prabhakar, T.R.: 1971, A singular integral equation with generalized

Mittag-Leffler function in the kernel, \emph {Yokohama Mathematical} 

\emph {Journal} {\bf 19}, 7-15.\\
Samko, S.G., Kilbas, A.A., and Marichev, O.I.: 1990, \emph {Fractional Integrals}

\emph {and Derivatives: Theory and Applications}, Gordon and Breach, 

New York.\\
Saxena, R.K., Mathai, A.M., and Haubold, H.J.: 2002, On fractional

kinetic equations, \emph {Astrophysics and Space Science} {\bf 282}, 281-287.\\
Saxena, R.K., Mathai, A.M., and Haubold, H.J.: 2004, On generalized

fractional kinetic equations, \emph {Physica A} {\bf 344}, 657-664.\\
Saxena, R.K., Mathai, A.M., and Haubold, H.J.: 2004a, Unified fractional

kinetic equation and a fractional diffusion equation, \emph {Astrophysics and}

\emph {Space Science} {\bf 290}, 299-310.\\
Saxena, R.K., Mathai, A.M., and Haubold, H.J.: 2004b, Astrophysical

thermonuclear functions for Boltzmann-Gibbs statistics and Tsallis

statistics, \emph {Physica A} {\bf 344}, 649-656.\\
Saxena, R.K., Mathai, A.M., and Haubold, H.J.: 2005, Fractional reaction-

diffusion equations, \emph {this volume}.\\
Smoller, J.: 1983, \emph {Shock Waves and Reaction-Diffusion Equations}, 

Springer-Verlag, New York-Heidelberg-Berlin.\\
Strier, D.E., Zanette, D.H., and Wio, H.S.: 1995, Wave fronts in a 

bistable reaction-diffusion system with density-dependent diffusivity, 

\emph {Physica A} {\bf 226}, 310.\\
Wilhelmsson, H. and Lazzaro, E.: 2001, \emph {Reaction-Diffusion Problems in} 

\emph {the Physics of Hot Plasmas}, Institute of Physics Publishing, 

Bristol and Philadelphia.\\
Wiman, A.: 1905, Ueber den Fundamentalsatz in der Theorie der

Functionen $E_\alpha(x)$, \emph {Acta Mathematica} {\bf 29}, 191-201.\\
\end{document}